\input amstex
\magnification=1200
\documentstyle{amsppt}
\NoRunningHeads
\NoBlackBoxes
\define\sothree{\operatorname{\frak s\frak o}(3,\Bbb C)}
\define\sltwo{\operatorname{\frak s\frak l}(2,\Bbb C)}
\define\RWsltwo{\operatorname{\Cal R\Cal W}(\sltwo)}
\define\Der{\operatorname{Der}}
\define\End{\operatorname{End}}
\define\CA{\Cal A}
\define\CF{\Cal F}
\define\Fg{\frak g}
\define\BC{\bold C}
\define\CS{\Cal S}
\define\CU{\Cal U}
\define\ga{\alpha}
\define\gb{\beta}
\define\gc{\gamma}
\define\gd{\delta}
\define\gep{\varepsilon}
\define\gz{\zeta}
\define\RWsltwoa{\operatorname{\Cal R\Cal W}(\sltwo;\alpha)}
\define\RWsltwoz{\operatorname{\Cal R\Cal W}(\sltwo;0)}
\define\CRWsltwo{\widehat{\operatorname{\Cal R\Cal W}}(\sltwo)}
\define\CRWsltwoa{\widehat{\operatorname{\Cal R\Cal W}}(\sltwo;\alpha)}
\define\slthree{\operatorname{\frak s\frak l}(3,\Bbb C)}
\define\sothreer{\operatorname{\frak s\frak o}(3)}
\define\RWsothreer{\operatorname{\Cal R\Cal W}(\sothreer)}
\define\tr{\operatorname{tr}}
\topmatter
\title Topics in hidden symmetries. IV.
\endtitle
\author Denis V. Juriev
\endauthor
\affil\rm
"Thalassa Aitheria" Research Center for Mathematical Physics and
Informatics,\linebreak
ul.Miklukho-Maklaya 20-180, Moscow 117437 Russia.\linebreak
E-mail: denis\@juriev.msk.ru\linebreak
\ \linebreak
November, 15, 1996\linebreak
E-version: q-alg/yymmxxx
\endaffil
\abstract This note being devoted to some aspects of the inverse problem of
representation theory explicates the links between researches on the Sklyanin
algebras and the author's (based on the noncommutative geometry) approach to
the setting free of hidden symmetries in terms of "the quantization of
constants".

Namely, the Racah--Wigner algebra for the Sklyanin algebra is constructed.
It may be considered as a result of the quantization of constants in the
Racah--Wigner algebra for the Lie algebra $\sltwo$. If the Sklyanin algebra
is interpreted as an algebra of anomalous spins then the Racah--Wigner algebra
for it may be regarded as an enlargement of the Sklyanin algebra by operators 
of the anomalous spin--spin interaction (of tensor type). 

The Racah--Wigner algebra for the Sklyanin algebra is an example of
the noncommutative weighted shift operator algebras (NWSO--algebras), which
generalize the mho--algebras introduced by the author earlier.
\endabstract
\endtopmatter
\document

This paper being a continuation of the previous three parts [1-3] as
illustrates the general ideology presented in the review [4] as explicates 
its new features. The subject is new links between a research activity on the
Sklyanin algebras [5,6] (see also [7]) and the author's (based on the 
noncommutative geometry) approach to the setting free of hidden symmetries 
[3:\S1] (see also [8,1]) in terms of "the quantization of constants" [4:\S3;3]. 
More concretely, we construct the Racah--Wigner algebra for the Sklyanin
algebra, which may be considered also as a result of the quantization of
constants in the Racah--Wigner algebra for the Lie algebra
$\sothree\simeq\sltwo$.

Formally, the Racah--Wigner algebra for the Sklyanin algebra is one of examples
of the noncommutative weighted shift operator algebras (NWSO--algebras), which
generalize the mho--algebras introduced by the author earlier [1:\S2] (see also 
[3:\S1.4]). 

Because the Sklyanin algebra is a very interesting object, deeply related to
many algebraic structures (see e.g.[9,10]), and the Racah--Wigner algebra for
the Lie algebra $\sltwo$ [8;11:\S2.2] (see also [3:\S\S1,2;2] and cf.[12:\S1.2])
is not less intriguing, their "hybrid" will be undoubtly of a certain value.

\head Topic Eight: The Racah--Wigner algebra for the Sklyanin algebra
\endhead

\subhead 1. Setting hidden symmetries free by the noncommutative Veronese 
mapping and the Racah--Wigner algebra $\RWsltwo$ for the Lie algebra $\sltwo$
\endsubhead

\definition{Definition 1 {\rm [8] (see also [1:\S1;4:\S1.1])}}

{\bf (A)} Let $\Fg$ be a lie algebra and $\CA$ be an associative algebra
such that $\Fg\subseteq\Der(\CA)$; a linear subspace $V$ of $\CA$ is called
{\it a space of hidden symmetries\/} iff (1) $V$ is a $\Fg$--submodule of
$\CA$, (2) the Weyl symmetrization defines a surjection $W:S^{\cdot}\mapsto\CA$
(the elements of $V$ are called {\it hidden symmetries with respect to\/}
$\Fg$). An associative algebra $\CF$ such that $\Fg\subseteq\Der(\CF)$ is
called {\it an algebra of the set free hidden symmetries\/} iff (1) $\CF$ is
generated by $V$, (2) there exists a $\Fg$--equivariant epimorphism of 
algebras $\CF\mapsto\CA$, (3) the weyl symmetrization defines an isomorphism
$S^{\cdot}(V)\mapsto\CF$ of $\Fg$--modules.

{\bf (B)} Let $V$ be a space of hidden symmetries and $\BC$ be some class of
associative Ore algebras; the hidden symmetries from $V$ are called 
$\BC$--regular iff the algebras of fractions $D(\CF)$ are isomorphic for
all those corresponding to $V$ algebras $\CF$ of the set free hidden
symmetries from the class $\BC$.

{\bf (C)} Let $V$ be a space of hidden symmetries in algebra $\CA$ with
respect to the Lie algebra $\Fg$; a subspace $V_0$ of $V$ is called {\it a
coordinate base\/} of $V$ iff (1) $V_0$ is a $\Fg$--submodule of $V$,
(2) the image of the Weyl symmetrization mapping $W_0:S^{\cdot}(V_0)\mapsto\CA$
contains $V$. Hidden symmetries from $V$ are called {\it of type\/} $(V_0,n)$
iff the image of $\oplus_{i\le n}S^i(V_0)$ under the weyl symmetrization
mapping $W_0$ coincides with $V$; in this case the mapping 
$\oplus_{i\le n}S^i(V_0)\mapsto\CF$, a composition of $W_0$ and the imbedding 
of $V$ into $\CF$, is called {\it the noncommutative Veronese mapping}. 

{\bf (D)} Let $\Fg$ be a Lie algebra, $V$ be a certain $\Fg$--module,
$\CA_s$ be a family of associative algberas, parametrized by $s\in\CS$ such 
that $\Fg\subseteq\Der(\CA_s)$, and $\pi_s:V\mapsto\CA_s$ be a family of
$\Fg$--equivariant imbeddings such that $\pi_s(V)$ is a space of hidden 
symmetries in $\CA_s$ with respect to $\Fg$ for a generic $s$ from $\CS$. An
associative algebra $\CF$ is called {\it an algebra of the
$\CA_{s,s\in\CS}$--universally set free hidden symmetries\/} iff $\CF$ is an
algebra of the set free hidden symmetries corresponding to 
$V\simeq\pi_s(V)$ for generic $\CA_s$ ($s\in\CS$). The hidden symmetries are 
called {\it $\CA_{s,s\in\CS}$--universally $\BC$--regular\/} iff the algebras
of fractions $D(\CF)$ are isomorphic for all algebras $\CF$ of the
$\CA_{s,s\in\CS}$--universally set free hidden symmetries from $\BC$.
\enddefinition

If $\Fg$ is a Lie algebra and $\CA$ is an associative algebra such that
$\Fg\subseteq\Der(\CA)$, and $V_0$ is a $\Fg$--submodule of $\CA$, whose
elements generate $\CA$ as an algebra, then in many interesting cases there
exists a space of hidden symmetries $V$ of type $(V_0,n)$ in $\CA$ for
a sufficiently large $n$.

\definition{Definition 2 {\rm [2] (cf. also [1:\S2;4:\S1.4])}}

{\bf (A)} Let $\Fg$ be a Lie algebra and $\pi$ be its (irreducible)
representation. {\it Mho--algebra\/} $\mho(\Fg,\pi)$ is an associative algebra
such that (1) $\CU(\Fg)$ is a subalgebra of $\mho(\Fg,\pi)$ and, hence, $\Fg$
naturally acts in $\mho(\Fg,\pi)$, (2) there is defined a $\Fg$--equivariant
embedding of $\pi$ into $\mho(\Fg,\pi)$, so $\pi$ may be considered as
a subspace of $\mho(\Fg,\pi)$, (3) the $\Fg$--equivariant embedding of $\pi$
into $\mho(\Fg,\pi)$ is extended to a $\Fg$--equivariant embedding of
$S^{\cdot}(\pi)$ into $\mho(\Fg,\pi)$, defined by the Weyl symmetrization,
and, therefore, $S^{\cdot}(\pi)$ may be considered as a subspace of
$\mho(\Fg,\pi)$; (4) $\Fg$--modules $\mho(\Fg,\pi)$ and $S^{\cdot}(\Fg)\otimes
S^{\cdot}(\pi)$ are isomorphic, here the isomorphism of subalgebra $\CU(\Fg)$
of the algebra $\mho(\Fg,\pi)$ and $S^{\cdot}(\Fg)$ as $\Fg$--modules is
used; (5) in an arbitrary basis $w_l$ in $\pi$ the commutator of two elements
of the basis in the algebra $\mho(\Fg,\pi)$ may be represented in the form
$[w_i,w_j]=f^k_{ij}w_k$, where the "noncommutative structural functions"
$f^k_{ij}$ are the elements of the algebra $\CU(\Fg)$.

{\bf (B)} Let $\Fg$ be a Lie algebra and $\pi$ be its (irreducible)
representation. {\it Affine mho--algebra\/} $\hat\mho(\Fg,\pi)$ is
an associative algebra such that the conditions (1)--(4) above hold and
(5') in an arbitrary basis $w_l$ in $\pi$ the commutator of two elements of the
basis in the algebra $\hat\mho(\Fg,\pi)$ may be represented in the form
$[w_i,w_j]=f^k_{ij}w_k+g_{ij}$, where the "noncommutative structural
functions" $f^k_{ij}$ and $g_{ij}$ are the elements of the algebra $\CU(\Fg)$.
\enddefinition

\proclaim{Theorem 1}

{\bf (A)} [8] (see also [1:\S1]) The tensor operators of type $\pi_1$ and
$\pi_2$ in the Verma module $V_h$ [$h$ is the weight] over the Lie algebra 
$\sltwo$ [$\pi_i$ is a finite--dimensional representation of $\sltwo$ of 
dimension $2i+1$] form a space of hidden symmetries of type $(\pi_1,2)$;
the quadratic (nonhomogeneous) algebras of the $\End(V_h)$--universally set
free hidden symmetries form an one--parametric family $\RWsltwoa$.

{\bf (B)} [1:\S2] (see also [4:\S1.4]) The Racah--Wigner algebra
$\RWsltwo=\mathbreak\RWsltwoz$ is a mho--algebra over $\sltwo$.

{\bf (C)} [8] (see also [1:\S1]) All the algebras of fractions $D(\RWsltwoa)$ 
are isomorphic [hence, the tensor operators of type $\pi_i$ ($i=1,2$) form a 
$\End(V_h)$--universally quadratic--regular scope of hidden symmetries].

{\bf (D)} [11:\S2.2] The Racah--Wigner algebra $\RWsltwo$ admits a quadratic
central extension $\CRWsltwo$, which may be considered also as a 
$\sltwo$--equivariant quadratic deformation of $\CU(\slthree)$, which is 
trivial on subalgebra $\CU(\sltwo)$.

{\bf (E)} [2] The central extension $\CRWsltwo$ of the Racah--Wigner 
algebra\linebreak $\RWsltwo$ is an affine mho--algebra over $\sltwo$.

{\bf (F)} [8] (see also [1:\S1]) The central extension $\CRWsltwo$ of
$\RWsltwo$ may be continued to the central extensions $\CRWsltwoa$ of
$\RWsltwoa$ in the class of quadratic (nonhomogeneous) algebras.
\endproclaim

The explicit formulas for the commutation relations in the algebras
$\RWsltwo$, $\CRWsltwo$, $\RWsltwoa$, $\CRWsltwoa$ were written in
[11:\S2.2;7] (see also [4:\S1.1]; a misprint in two commutators in [11:\S2.2]
was corrected in [8]).

Let's consider a "compact" real form of the Racah--Wigner algebra
$\RWsltwo$, which will be denoted by $\RWsothreer$. The compactness means that
the algebra $\RWsothreer$ is a mho--algebra over $\sothreer$.

\proclaim{Theorem 2} The Racah--Wigner algebra $\RWsothreer$ is generated
by eight generators $s_i$, $t_{ij}=t_{ji}$ ($i,j=1,2,3$, 
$t_{11}+t_{22}+t_{33}=0$) with commutation relations
$$\aligned
[s_i,s_j]&=\epsilon_{ijk}s_k,\quad [s_i,t_{jk}]=\epsilon_{ijl}t_{lk}+
\epsilon_{ikl}t_{lk},\\
[t_{ij},t_{kl}]&=\epsilon_{ikm}(s_jt_{ml}+s_lt_{jm})+
\epsilon_{ilm}(s_jt_{mk}+s_kt_{jm})\\
&+\epsilon_{jkm}(s_it_{ml}+s_lt_{im})+\epsilon_{jlm}(s_it_{mk}+s_kt_{im}),
\endaligned$$
where $\epsilon_{ijk}$ is a totally antisymmetric tensor.
\endproclaim

\remark{Remark 1} Note that one may use the left, the right or the Weyl
ordering of the r.h.s. of commutation relations between $t_{ij}$, there
are no differences between them.
\endremark

\remark{Remark 2} The relations of the Theorem 2 between $s_i$ and $t_{jk}$
may be considered as classical, they define correct Poisson brackets.
\endremark

\subhead 2. The quantization of constants and the Sklyanin algebra
\endsubhead

Under {\it the quantization of constants\/} one means the change of certain
numerical constant in the commutation relations to the matrix (operator)
quantities. An attempt of rigorous definition of the procedure was made in
[3:\S1].

The Lie algebra $\sothreer$ is realized as a Lie algebra of all
skew--symmetric $3\times 3$--matrices. Let's fix a basis
$$e_1=\left(\matrix 0 & 1 & 0 \\ -1 & 0 & 0 \\ 0 & 0 & 0 \endmatrix\right),\quad
e_2=\left(\matrix 0 & 0 & 0 \\ 0 & 0 & 1 \\ 0 & -1 & 0 \endmatrix\right),\quad
e_3=\left(\matrix 0 & 0 & 1 \\ 0 & 0 & 0 \\ -1 & 0 & 0 \endmatrix\right)$$
in the algebra. The commutation relations in such basis are of the form
$$[e_i,e_j]=\epsilon_{ijk}e_k,$$
where $\epsilon_{ijk}$ is a totally antisymmetric tensor. 

\proclaim{Theorem 3 {\rm [4:\S3.1]}} For any triple of independent
skew--symmetric $3\times 3$--matrices $S_i$ ($i=1,2,3$) there exists only
one symmetric $3\times 3$--matrix $Q$ such that
$$[S_i,S_j]=\epsilon_{ijk}(QS_k+S_kQ),$$
where $\epsilon_{ijk}$ is a totally antisymmetric tensor. Moreover,
$$[Q,S_i]=F_i(S_1,S_2,S_3),$$
where $F_i$ are uniquely defined quadratic forms of $S_1$, $S_2$, $S_3$.
If matrices $S_i$ are written in the form
$$S_1=\left(\matrix 0 & \ga & 0 \\ -\ga & 0 & 0 \\ 0 & 0 & 0 \endmatrix\right),\quad
S_2=\left(\matrix 0 & \gb & \gc \\ -\gb & 0 & 0 \\ -\gc & 0 & 0 \endmatrix\right),\quad
S_3=\left(\matrix 0 & \gd & \gep \\ -\gd & 0 & \gz \\ -\gep & -\gz & 0 \endmatrix\right),$$
then
$$Q=\left(\matrix u & x & y \\ x & v & z \\ y & z & w \endmatrix\right),$$
where
$$\align
x=\frac{\ga^2\gep+\gb^2\gep-\gb\gc\gd}{\ga\gc},&\quad
y=\frac{\gb\gep-\gc\gd}{\ga},\quad
z=\frac{\gb\gz}{\ga},\\
u+v=-\frac{\gc\gz}{\ga},\quad
u+w=-\tfrac{\gz(\ga^2+\gb^2)}{\ga\gc},&\quad
v+w=\tfrac{2\gb\gc\gd\gep-\ga^2\gep^2-\gb^2\gep^2-\ga^2\gc^2-\gep^2\gc^2}{\ga\gc\gz}.
\endalign$$
Matrices $S_i$, $Q$ realize a representation of a quadratic algebra, which
naturally coincides with the Sklyanin algebra [5,6] if $\tr(S_iS_j)=0$ for 
$i\ne j$. In the case of the Sklyanin algebra 
$F_i(S_1,S_2,S_3)=\epsilon_{ijk}J_{jk} 
(S_jS_k+S_kS_j)$ ($J_{12}+J_{23}+J_{31}+J_{12}J_{23}J_{31}=0$). 
\endproclaim

\remark{Remark 3} $Q$ is nondegenerate. $\tr(S_iQS_j+S_jQS_i)=0$ if
$i\ne j$ (i.e. $S_i$ are $Q$--orthogonal).
\endremark

\remark{Remark 4} (A) Let $Q$ be an arbitrary nondegenerate symmetric matrix.
Put $p(Q)A=QA+AQ$. If $S_1$ and $S_2$ are two skew--symmetric $Q$--orthogonal 
matrices, then $S_3=p(Q)^{-1}([S_1,S_2])$ is $Q$--orthogonal to both
$S_1$ and $S_2$. (B) If $Q$ is a nondegenerate symmetric $3\times 3$--matrix
and $S_1$, $S_2$ be two $Q$--orthogonal skew--symmetric $3\times 3$--matrices,
then there exist three multiplies $\lambda_1$, $\lambda_2$ and $\lambda_3$
such that $\lambda_1S_1$, $\lambda_2S_2$, $\lambda_3p(Q)^{-1}([S_1,S_2])$ and
$Q$ generate the algebra of the theorem 3.
\endremark

\remark{Remark 5} The algebra of the theorem 3 is generally isomorphic to the
Sklyanin algebra (this is a manifestation of the fact that any symmetric
matrix may be written in a diagonal form).
\endremark

\subhead 3. The Racah--Wigner algebra for the Sklyanin algebra
\endsubhead

\definition{Definition 3} 

{\bf (A)} Let $\CA$ be an associative algebra generated
by the subspace $V$ so that the Weyl symmetrization mapping 
$W:S^{\cdot}(V)\mapsto\CA$ is a bijection. {\it A noncommutative weighted
shift operator algebra (NWSO--algebra)\/} over $\CA$ is an associative algebra
$\CF$ generated by $\hat V$ such that (1) the Weyl symmetrization mapping
$\hat W:\S^{\cdot}(\hat V)\mapsto\CF$ is a bijection, (2) the subspace
$\hat V$ may be decomposed into the direct sum $V\oplus U$, where the
elements of $V$ generate a subalgebra isomorphic to $\CA$, (3)
there exist a mapping $R:U\otimes V\mapsto\CA\otimes U$ ($R(t\otimes 
s)=\sum_i \varrho^{(i)}_1(s,t)\otimes \varrho^{(i)}_2(s,t)$) such that
$st=\sum_i\varrho^{(i)}_1(s,t)\varrho^{(i)}_2(s,t)$, (4) there exists a
mapping $T:\Lambda^2(U)\mapsto\CA\otimes U$ ($T(t_1\wedge t_2)=
\sum_i\tau^{(i)}_1(t_1,t_2)\otimes\tau^{(i)}_2(t_1,t_2)$) such that
$[t_1,t_2]=\sum_i\tau^{(i)}_1(t_1,t_2)\tau^{(i)}_2(t_1,t_2)$.

{\bf (B)} Let $\CA$, $V$ be the same data as in (A). {\it An affine 
noncommutative weighted shift operator algebra (affine NWSO--algebra)\/}
over $\CA$ is an associative algebra $\CF$ generated by $\hat V$ such that
the conditions (1)-(3) above hold, and (4') there exist two mappings
$T:\Lambda^2(U)\mapsto\CA\otimes U$ ($T(t_1\wedge t_2)=
\sum_i\tau^{(i)}_1(t_1,t_2)\otimes\tau^{(i)}_2(t_1,t_2)$) 
and $S:\Lambda^2(U)\mapsto\CA$ ($S(t_1\wedge t_2)=\sigma(t_1,t_2)$) such that
$[t_1,t_2]=\sum_i\tau^{(i)}_1(t_1,t_2)\tau^{(i)}_2(t_1,t_2)+\sigma(t_1,t_2)$.
\enddefinition

Definition 3 define objects, which are generalizations of mho--algebras and
affine mho--algebras of Definition 2. Note that Definitions 1-3 allows a
slight generalization in terms of the Poincar\'e--Birkhoff--Witt algebras
[10:App.1]. 

\proclaim{Lemma 1} Let the Sklyanin algebra be realized as in the theorem 3.
Then there exist six symmetric $Q$--orthogonal matrices $T_{ij}=T_{ji}$
($i,j=1,2,3$) such that $[S_i,T_{jk}]=\epsilon_{ijl}P(Q)T_{lk}+
\epsilon_{ikl}P(Q)T_{jl}$ ($\epsilon_{ijk}$ is a totally antisymmetric
tensor), moreover, $Q$ is proportional to $T_{11}+T_{22}+T_{33}$. The matrices
$T_{ij}$ ($i\ne j$) and $T_{ii}-T_{jj}$ are uniquely defined up to a common
multiple.
\endproclaim

\proclaim{Lemma 2} Let the Sklyanin algebra be generated by the elements $S_i$
($i=1,2,3$), $Q$. Let $T_{ij}$ be the formal variables, the commutation 
relations of which with $S_i$ are of the form prescribed in Lemma 1. Then there 
exist a uniquely defined mapping $\Xi$ 
[ $\Xi(T_{ij})=\sum_{klm}\xi^{klm}_{ij}S_kT_{lm}+T_{lm}S_k$ ] defined on $T_{ij}$
($i\ne j$) and $T_{ii}-T_{jj}$ such that if one puts $[Q,T_{ij}]=\Xi(T_{ij})$ 
then the Jacobi identities for triples $Q$, $S_i$, $T_{jk}$ and $S_i$, $S_j$, 
$T_{kl}$ will formally hold.
\endproclaim

\proclaim{Theorem 4} Let the Sklaynin algebra be generated by the elements $S_i$
($i=1,2,3$), $Q$ and realized as in Theorem 3. Let $T_{ij}$ be the formal 
variables realized as in Lemma 2 with the specified there commutation relations.
Then there exist a unique NWSO--algebra over the Sklyanin algebra with
the subspace $U$ generated by $T_{ij}$ ($i\ne j$) and $T_{ii}-T_{jj}$
compatible with the realizations of $S_i$, $Q$ and $T_{jk}$ and with the
mapping $T$ being linear on $S_i$.
\endproclaim

\remark{Remark 6} The algebra of the Theorem 4 will be called the Racah--Wigner
algebra for the Sklaynin algebra. It may be considered as a result of the
quantization of constants in the Racah--Wigner algebra $\RWsothreer$. On the
other hand, its construction is similar to the procedure of the setting free of
hidden symmetries above and may be regarded as its quantum analogue. 
\endremark

\remark{Remark 7} The construction of the Racah--Wigner algebra for the 
Sklyanin algebra has its classical counterpart in view of the remark 2.
\endremark

\remark{Questions}

(1) Has the Racah--Wigner algebra for the Sklyanin algebra a central extension
and one--parametric deformation, which are results of the quantization of
constants in objects of the Theorem 1?

(2) Is the construction of [1:\S2] (see also [4:\S1.4]) applicable to the
Sklyanin algebra?

(3) Does the construction of the Racah--Wigner algebra for the Sklyanin algebra
admit a superanalogue (see [11:\S2.3] on a superanalogue of $\RWsltwo$)?

(4) What are the generalized Mickelsson and Zhelobenko algebras of the
constructed Racah--Wigner algebra for the Sklyanin algebra (cf.[2:\S1])?

(5) The construction above is formally inapplicable to the degenerate
Sklyanin algebra [6,7] (see also [10:App.1]. It is interesting to formulate
a generalization or to adapt the construction to this specific case.

(6) Does the Sklyanin algebra admit an analogue of the construction of
[1:\S1] (see also [4:\S1.3]) of the setting free of hidden symmetries by
the noncommutative semicubic mapping?

(7) What is a relation between the Racah--Wigner algebra for the Sklyanin
algebra and its classical counterpart?

(8) The Racah--Wigner algebra for $\sltwo$ may be constructed in lines of
[4:\S1.6] as generated by the operators of averaged spin and averaged
spin--spin interaction. Such procedure was called by {\it a resynthesis of 
algebraic structure}. The startpoint of the procedure is $\CU(\sltwo)$ to 
obtain $\RWsltwo$. Is the resynthesis applicable to the Sklyanin algebra to 
receive the Racah--Wigner algebra for it? A supposition of applicability is 
supported by arguments based on a physical intuition, which allows to identify 
the Racah--Wigner algebra for the Sklyanin algebra with an enlargement of the
least by operators of the anomalous spin-spin interaction (of tensor type 
[13:Ch.XVI,\S116]) if the Sklyanin algebra is considered as an algebra of 
anomalous spins.
\endremark

\Refs
\roster
\item" [1]" Juriev D., Topics in hidden symmetries. I. E-print: hep-th/9405050.
\item" [2]" Juriev D., Topics in hidden symmetries. II. Report RCMPI-96/06
(October 1996) [e-version: q-alg/9610026 (1996)].
\item" [3]" Juriev D., Topics in hidden symmetries. III. Report RCMPI-96/06$^+$
(November 1996) [e-version: q-alg/9611003 (1996)].
\item" [4]" Juriev D., An excursus into the inverse problem of representation
theory [in Russian]. Report RCMPI-95/04 (August 1995) 
[e-version: mp\_arc/96-477].
\item" [5]" Sklyanin E.K., On algebraic structures related to the Yang--Baxter
equation [in Russian]. Funkts. analiz i ego prilozh. 16(4) (1982) 27-34; 17(4)
(1983) 34-48.
\item" [6]" Sklyanin E.K., On an algebra generated by the quadratic relations
[in Russian]. Uspekhi Matem. Nauk 40(2) (1985) 214.
\item" [7]" Tarasov V.O., Takhtadzhan L.A., Faddeev L.D., Local hamiltonians
for integrable quantum models on a lattice [in Russian]. Teor.Matem.Fiz.
57(2) (1983) 163-181.
\item" [8]" Juriev D., Setting hidden symmetries free by the noncommutative
Veronese mapping. J.Math.Phys. 35(9) (1994) 5021-5024.
\item" [9]" Reshetikhin N.Yu., Takhtadzhan L.A., Faddeev L.D., Quantization of
Lie groups and Lie algebras [in Russian]. Algebra i analiz 1(1) (1989) 178-206
[English transl.: St.Petersburg Math.J. 1 (1990) 193-225].
\item"[10]" Karasev M.V., Maslov V.P., Nonlinear Poisson brackets. Geometry and
applications [in Russian]. Moscow, Nauka, 1991 [English transl. by the 
Amer.Math.Soc., R.I., 1993].
\item"[11]" Juriev D., Complex projective geometry and quantum projective field
theory [in Russian]. Teor.Matem.Fiz. 101(3) (1994) 331-348 [English transl.:
Theor.Math.Phys. 101 (1994) 1387-1403].
\item"[12]" Juriev D., Infinite dimensional geometry and quantum field theory
of strings. III. Infinite dimensional W--geometry of a second quantized free
string. J.Geom.Phys. 16 (1995) 275-300.
\item"[13]" Landau L.D., Lifschitz E.M., A course of theoretical physics. III.
Quantum mechanics [in Russian]. Moscow, Fizmatlit, 1963.
\endroster
\endRefs
\enddocument